\documentclass[11pt]{amsart}

\usepackage[english]{babel}

\usepackage{amssymb,latexsym,amsmath,amsthm,graphicx}
\newtheorem{thm}{Theorem}
\newtheorem{lemma}[thm]{Lemma}

\newtheorem*{claim}{Claim}

\theoremstyle{remark} 
\newtheorem{remark}[]{Remark}

\usepackage{xcolor}

\newcommand{\NN}{\mathbb{N}}

\newcommand{\RR}{\mathbb{R}}
\newcommand{\EE}{\mathbb{E}}
\newcommand{\PP}{\mathbb{P}}
\newcommand{\CC}{\mathbb{C}}
\newcommand{\f}{f}
\newcommand{\HH}{H_{[\alpha d,d]}}
\newcommand{\B}{B_{x}}
\newcommand{\RP}{\mathbb{R}\textrm{P}}
\newcommand{\CP}{\mathbb{C}\textrm{P}}

\begin{document}
\author{Antonio Lerario and Erik Lundberg}
\title{Statistics on Hilbert's Sixteenth Problem}
\maketitle
\begin{abstract}
We study the statistics of the number of connected components and the volume of a random real algebraic hypersurface in $\RP^n$ defined by 
a homogeneous polynomial $f$ of degree $d$ in the \emph{real Fubini-Study} ensemble. We prove that for the expectation of the number of connected components:
\begin{equation} \label{bettiabstract}\EE b_0(Z_{\RP^n}(f))=\Theta(d^n),\end{equation}
the asymptotic being in $d$ for $n$ fixed.\\
We do not restrict ourselves to the random homogeneous case and we consider more generally random polynomials belonging to a window of eigenspaces of the Laplacian on the sphere $S^n,$ proving that the same asymptotic holds. As for the volume properties we prove that:
\begin{equation}\label{volumeabstract}\EE \textrm{Vol}(Z_{\RP^n}(f))=\Theta(d).\end{equation} 
Both equations (\ref{bettiabstract}) and (\ref{volumeabstract}) exhibit expectation of maximal order in light of  Milnor's bound $b_0(Z_{\RP^n}(f))\leq O(d^n)$ and the bound $\textrm{Vol}(Z_{\RP^n}(f))\leq O(d).$

\end{abstract}
\section{Introduction}

Hilbert's sixteenth problem, in its general form, asks for the study of the maximal number and the possible arrangements of the components of a nonsingular real algebraic hypersurface of degree $d$ in $\RP^n$. 
It turns out that even the case of  curves in $\RP^2$ is an extremely subtle problem and essentially nothing is known for $n\geq 3$. 
The analogous question for a hypersurface in $\CP^n$ is trivial, as for a \emph{generic} choice of the polynomial defining it, its topology is determined 
(this is the statement that the cohomology of a complete intersection in the complex projective space is determined by its multidegree and complete intersections are generic).

Thus, if from one hand no technique exists to study the arrangement of the components of a real algebraic hypersurface, on the other hand one still would like to formulate ``typical" statements: 
a way for doing that is to replace the word \emph{generic} with \emph{random}.

Let us start by considering the case of hypersurfaces in $\RP^1$, i.e. zeros of ordinary univariate polynomials. 
In fact here we consider the \emph{homogenization} of a univariate polynomial together with its \emph{projective} roots rather than the classical polynomial itself, 
but it is clear that the two objects produce the same statistic.
Over the complex numbers statements like:
$$\emph{``a generic complex polynomial of degree $d$ has $d$ roots"}$$ make perfect sense. 
On the other hand the only property we can say to be generic for the real roots of a real polynomial is the one of being distinct.

Proposing a random viewpoint raises the obvious question: \emph{what is meant by random?} 
We first look to the seminal work \cite{Kac} of Kac, who proved that a random univariate polynomial of degree $d$ whose coefficients are i.i.n. has asymptotically $2/\pi\log d$ real roots.  
Even if it might seem natural, the distribution Kac considered is not invariant by the action of projective transformations over the extended real line. 
To produce such an invariant distribution we proceed as follows. 
We first \emph{extend} our space of polynomials by homogenization; 
this leads us to consider \emph{projective} zeros rather than just on the real line. 
Thus, we consider the vector space $W_{1,d}$ of real homogeneous polynomials of degree $d$ in two variables.
We endow this vector space with a Euclidean structure by considering the $L^2$-norm of the restriction of these polynomials to the circle $S^1$. In an equivalent way we can define for $f$ and $g$ in $W_{1,d}$ 
their scalar product to be:
$$\langle f, g\rangle=\int_{\RR^2} f(x)g(x)e^{-\|x\|^2}dx,$$
where $x$ is the variable in $\RR^2$. Notice that since $f$ and $g$ are homogeneous of the same degree, the previous integral equals the integral of $fg$ restricted to the circle $S^1$ (up to a constant factor).
Using this scalar product we can define a \emph{centered Gaussian} probability distribution on $W_{1,d}$ by letting for a measurable subset $A$:
$$\textrm{probability of $A$}=\frac{1}{\nu_{1,d}}\int_A e^{-\frac{\|f\|^2}{2}}df,$$
where $\nu_{1,d}=(2\pi)^{(d+1)/2}$ is the constant that makes the previous a probability distribution and $f$ is treated as the integration variable that runs over $W_{1,d}.$

For us a \emph{random polynomial} will simply be a polynomial in $W_{1,d}$ with the above distribution.

Having thus replaced \emph{generic} with \emph{random}, we can make a positive statement in the real setting:
$$\emph{``a real random polynomial of degree $d$ has $\sqrt{\frac{d(d+2)}{3}}$ roots".}$$
We provide now another equivalent definition of a random polynomial.
We fix a basis $\{f_1, \ldots, f_{d+1}\}$ for $W_{1,d}$ which is orthonormal for the above scalar product and we consider a random function of the form:
$$f=\sum_{i=1}^{d+1}\xi_i f_i\quad\textrm{with}\quad \xi_i\sim N(0,1).$$
What we get in this way is again a random polynomial in the above sense.
One way to produce such a basis is by considering \emph{trigonometric polynomials}, i.e. functions on the circle satisfying:
$$\Delta_{S^1}f=-l^2 f,\quad l\in \NN.$$ 
We denote by $H_{1,l}$ the vector space of all the solutions of the previous equation: these must be linear combinations of $\sin (l \theta)$ and $\cos (l\theta)$. 
Equivalently we can extend such functions by homogeneity to all $\RR^2$ and we get harmonic polynomials of degree $l$. 
What is interesting for us is that in this fashion we can write the vector space $W_{1,d}$ as:
$$W_{1,d}=\bigoplus _{d-l \in 2\NN}\|x\|^{d-l}H_{1,l},$$ 
where again $x$ is the variable in $\RR^2$ and the direct sum turns out to be orthogonal with respect to the above scalar product.

In particular the above orthonormal basis can be obtained by collecting together the orthonormal bases coming from each $H_{1,l}$, $d-l\in 2\NN$.

If we repeat the same construction for the scalar product 
$$\langle f, g\rangle_{\CC} = \int_{\CC^2}f(z)\overline{g(z)}e^{-\|z\|^2}dz,$$
then the corresponding real random polynomials are called \emph{Kostlan} distributed (or \emph{Bombieri-Weyl}). The vector space is the same, but in this second case the Euclidean structure comes by restricting to the real polynomials the Hermitian structure given on the complex ones. The expectation of the number of real roots of a real Kostlan polynomial of degree $d$ is $\sqrt{d}$. 
Already from this comparison we see that the structure of the zeros of a random polynomial is richer than the one of a Kostlan, which in turn is richer than for the naive ensemble.
\begin{remark}The two numbers $d(d+2)/3$ and $d$ are called the \emph{parameters} of these distributions. More generally each Gaussian distribution on $W_{1,d}$ that is invariant by projective transformation has its parameter $\delta$ which is defined to be the square of the expected number of zeros it produces.\end{remark}
Our interest is in what happens when we increase the number of variables. We let $W_{n,d}$ be the vector space of homogeneous polynomials of degree $d$ in $n+1$ variables,
and consider their zero sets as hypersurfaces in $\RP^n$.
If in the univariate case as a ``measure" of our set we could use the number of its points, 
here we have at least two possible approaches to generalize this idea.

The first one is by considering the number of points as a zero-dimensional volume. In this case an upper bound for the ``size" of our set is provided by\footnote{Here the notation $f(d)=\Theta(g(d))$ means that there exist two constants $c_1, c_2$ such that for $d$ large enough $c_1g(d)\leq f(d)\leq c_2 g(d).$}:
 \begin{equation}\label{introvolume}\textrm{Vol}(Z_{\RP^n}(f))\leq \Theta( d).\end{equation}
\begin{remark}The previous formula comes from an application of the Integral Geometry formula and Wirtinger's inequality (see \cite{Moncet}). In fact if $f$ is regular the
previous volume is bounded by $d\cdot\textrm{Vol}(\RP^{n-1})/\textrm{Vol}(\CP^{n-1})$ and since $d=\textrm{Vol}(Z_{\CP^n}(f))$ this also gives a nice way to compare the real zero
locus to the complex one (as we did in the univariate case using the Fundamental Theorem of Algebra).\end{remark}
 
The second approach is more subtle and regards the number of points as a topological property (rather than a metric one), namely as the number of connected components. 
It is interesting to notice that there is a quantitative bound also for this number, usually referred to as Milnor's bound (see \cite{Milnor}):
\begin{equation}\label{Milnor}b_0(Z_{\RP^n}(f))\leq \Theta(d^{n}).\end{equation}

\begin{remark}The original paper \cite{Milnor} of Milnor does not contain the bound in this form, since a $d^{n+1}$ term appears. We will give a proof of (\ref{Milnor}) in the Appendix for the smooth case (which is the case we are interested in, since the set of nonregular polynomials has \emph{measure zero} in $W_{n,d}$), though the same bound holds true even for singular hypersurfaces.
\end{remark}
As we already pointed out, for a generic $f$ in $W_{n,d}$ both the volume and the Betti numbers are constant for its complex zero locus. This is certainly not the case for the real zero locus, as for example it might even be empty. To study the expectation of these quantities over the reals we make $W_{n,d}$ into a Gaussian space. We first define the scalar product of two elements $f$ and $g$ in it to be:
$$\langle f, g\rangle=\int_{\RR^{n+1}}f(x)g(x)e^{-\|x\|^2}dx.$$
Using this scalar product $W_{n,d}$ is turned into a probability space by defining for each of its measurable subsets $A$:
$$\textrm{probability of $A$}=\frac{1}{\nu_{n,d}}\int_A e^{-\frac{\|f\|^2}{2}}df,$$
where now $\nu_{n,d}=(2\pi)^{\textrm{dim}(W_{n,d})/2}.$ 

Following \cite{Sarnak2011}, we call $W_{n,d}$ with the resulting distribution of probability the \emph{real Fubini-Study ensemble}.

The expectation of the volume of the zero locus of a random polynomial of degree $d$ from the real Fubini-Study ensemble was computed by P. B\"urgisser in \cite{Burgisser2007}:
\begin{equation}\label{volumeburg}\EE \textrm{Vol}(Z_{\RP^n}(f))=\sqrt{\frac{d(d+2)}{3}}\textrm{Vol}(\RP^{n-1}).\end{equation}
In fact in the paper \cite{Burgisser2007} this result follows from the more striking computation of the expected curvature polynomial of $Z_{\RP^n}(f)$ for any projectively
invariant distribution. Such a result is formulated in terms of the above mentioned parameter $\delta$, which in this case is the square of the average number of real zeros of a restriction to an $\RP^1$.

In particular, if we define the \emph{Kostlan} distribution on $W_{n,d}$ by means of the scalar product\footnote{The corresponding norm is usually called the \emph{Bombieri norm}.} $\int_{\CC^{n+1}}f(z)\overline{g(z)}e^{-\|z\|^2}dz $, 
the corresponding parameter is $d$ and the result we get for the expectation of the volume is 
$$\EE \textrm{Vol}(Z_{\RP^n}(f))=\sqrt{d}\textrm{Vol}(\RP^{n-1})$$
(the general result for a distribution with parameter $\delta$ is $\delta^{1/2}\textrm{Vol}(\RP^{n-1})$).

In the case $n=2$ and restricting the distribution of probability to $H_{2,d}$, i.e. considering only \emph{harmonic} polynomials rather than homogeneous ones, 
the statistics of the number of connected components was studied by Nazarov and Sodin. 
In the breakthrough paper \cite{NazarovSodin2009} the authors prove that there exists a constant $c>0$ such that for a random spherical harmonic $f$ of degree $d$:
$$\lim_{d\to \infty}\frac{\EE b_0(Z_{S^2}(f))}{d^2}=c.$$
Moreover, they showed that $b_0$ concentrates exponentially about its mean.

In the visionary handwritten letter \cite{Sarnak2011} P. Sarnak claims that using the same ideas as in \cite{NazarovSodin2009} it is possible to extend this result to the case of  real homogeneous polynomials in the real Fubini-Study ensemble. The interesting part of these statements (besides the concentration) is the lower bound estimate: 
an upper bound of the same order is given by Milnor's one (or in the spherical harmonic case by Courant's nodal domain theorem). The statistics of the number of connected components of a Kostlan distributed curve was studied by Gayet and Welschinger: 
in the paper \cite{GayetWelschinger2012} the authors prove that such a real curve of degree $d$ has on average $O(d)$ components. 
The same authors prove in \cite{GayetWelschinger3} that the expectation of \emph{each} Betti number of a Kostlan hypersurface in $\RP^n$ is bounded by $O(d^{n/2})$, 
generalizing in some sense the univariate case.

On a different perspective the first author in \cite{Lerario2012}, using techniques from Random Matrix Theory, has studied intersection of one or two Kostlan distributed quadrics, the asymptotic being this time in the number $n$ of variables; for such random algebraic sets the expectation of the sum of the Betti numbers behaves asymptotically as $n$.

As we did before, it is possible to give an equivalent definition of the real Fubini-Study structure on $W_{n,d}$ using an $L^2(S^n)$-orthonormal basis. In fact we can again decompose our vector space into the \emph{orthogonal} sum:
$$W_{n,d}=\bigoplus _{d-l \in 2\NN}\|x\|^{d-l}H_{n,l},$$
where now $H_{n,l}$ is the vector space of spherical harmonics of degree $l$, i.e. those functions on the sphere $S^n$ satisfying:
$$\Delta_{S^n}f=-l(l+n-1)f.$$

Inspired by the ideas of \cite{NazarovSodin2009} and \cite{Sarnak2011} we study a  more general class of random polynomials that contains as special cases both the homogeneous and the harmonic ones. 
To define such a class, for every $\alpha \in [0,1]$ and $d>0$ let us consider the following index set:
$$[\alpha d, d]=\{\textrm{natural numbers $l$ such that $[\alpha d]\leq l\leq d$ and $d-l$ is even}\}.$$
With this definition the $\alpha$-\emph{window} of eigenspaces we are going to consider is the following orthogonal direct sum:
$$\HH=\bigoplus _{l\in [\alpha d, d]}H_{n,l},\quad \alpha\in [0,1].$$
We endow $H_{[\alpha d, d]}$ with the probability distribution induced by the real Fubini-Study ensemble. 

\begin{remark}There is no unified convention for the names of these ensembles. 
In fact the name \emph{real} Fubini-Study should be distunguished from \emph{complex} Fubini-Study
which is another name for what we called Kostlan distributed.\\
In this paper we will only deal with the \emph{real Fubini-Study} ensemble and its restriction to $H_{[\alpha d, d]}.$ For this reason, since no confusion can arise, a \emph{random polynomial} for us will simply mean a polynomial from this ensemble.\end{remark}

In particular for $\alpha=0$ we get random homogeneous polynomials of degree $d$ and for $\alpha=1$ we get  random spherical harmonics of degree $d$: 
$$H_{[0,d]}= W_{n,d}  \quad \textrm{and}\quad H_{[d,d]}=H_{n,d}.$$ 

\begin{remark}
In fact we do not need to restrict to the case the $\alpha$-window contains only the eigenspaces $H_{n,l}$ with $d-l$ even. 
Our choice is motivated by the fact that we want to avoid redundancy in order to get random polynomials for $\alpha=0$; one can consider as well the full window.\end{remark}
We denote the dimension of $H_{n,l}$ by $d(n,l)$ and the one of $\HH$ by $D_\alpha$; notice that in the case $\alpha\neq 1$ we have $D_\alpha=\Theta(d^n)$. 
If we consider an orthonormal basis $\{Y_l^i\}_{i=1}^{d(n,l)}$ for $H_{n,l}$, we can decompose a random $\f$ in $\HH$ into a sum of random spherical harmonics of different degrees:
$$\f=\sum_{l\in[\alpha d, d]}\sum_{i=1}^{d(n,l)}\xi_{l}^iY_l^i\quad \textrm{with}\quad \xi_l^i\sim N(0, D_\alpha^{-1}).$$
The introduction of the scaling coefficient $D_\alpha^{-1}$ is absolutely not necessary,
but it has the advantage that $\EE \|f\|^2_{L^2(S^n)}=1$ which simplifies the notation in the proofs.

As for the expectation of metric properties of the zero locus of a random $f$ in $\HH$, we prove the following theorem.
\begin{thm}\label{volume}Let $\alpha \in [0,1]$ and $f$ be a random polynomial in $H_{[\alpha d, d]}$. Then:
$$\EE \emph{\textrm{Vol}}(Z_{S^n}(f))=\Theta (d).$$
\end{thm}
In particular taking expectation on both sides of inequality (\ref{introvolume}) turns it into an equality (the volume on the projective space and the one on the sphere being related by a constant that depends only on $n$).

We notice that the statement for $\alpha=1$ was proved by B\'erard in \cite{Berard} 
and for $\alpha=0$ by B\"urgisser \cite{Burgisser2007} as stated above in (\ref{volumeburg}). 
In fact it is possible to write the exact value of this expectation using an analog of what we called the parameter of the distribution, defined in this case as $\delta=\frac{1}{nD_\alpha}\sum_{l\in [\alpha d, d]}l(l+n-1)d(n,l)$; the expectation of the volume is then $\delta^{1/2}\textrm{Vol}(S^{n-1}).$ The same statement holds for the projective zero locus, as it is double covered by the spherical one.

The statistics of the number of connected components is more subtle and requires more work, drawing on the ideas of \cite{NazarovSodin2009}. The statement we will prove is the following. \begin{thm}\label{thm:main}
Let $\alpha\in [0,1]$, and $f$ be a random polynomial in $\HH$. Then:
$$\EE b_0(Z_{S^n}(f))=\Theta(d^n).$$\end{thm}
In analogy with the case of the volume, this theorem says that if we take expectation on both sides of Milnor's bound (\ref{Milnor}) we turn it into an equality (the number of connected components on the projective space and those on the sphere being related by a factor of \emph{two}).

Notice that the case $n=2$ and $\alpha=1$ gives the spherical harmonic case proved by Nazarov and Sodin in \cite{NazarovSodin2009}. The case $n=2$ and $\alpha=0$ gives Sarnak's claim \cite{Sarnak2011} on curves from the real Fubini-Study ensemble.

\begin{remark}After this paper was written, the authors have been informed by F. Nazarov and M. Sodin that they have an alternative approach to this problem, that works in more general settings; a summary of their ideas is available at \cite{slides}.

The authors have also been informed by P. Sarnak that together with I. Wigman they are able to provide asymptotics for the expected number of ovals of a random plane curve for the \emph{naive} model (the defining coefficients for the monomials are independent normals with mean zero and variance one), the Kostlan model and the real Fubini-Study one. 
\end{remark}

\subsection{A random version of Hilbert's sixteenth problem}

It is interesting to discuss the case $\alpha=0$ and $n=2$, i.e. random real algebraic curves. The first part of Hilbert's sixteenth problem asks to study the configuration of the components of a real algebraic curve of degree $d$ in the plane $\RP^2.$ The nullhomotopic components of such curves are called \emph{ovals}; each one of them separates $\RP^2$ in two components: one of them is homeomorphic to a disk (and is called the \emph{interior} of the oval) and the other to a Moebius band. In the case $d$ is even all components are ovals, but if $d$ is odd there is a component that is not an oval (it is homotopic to a projective line and generates the fundamental group of $\RP^2$). Harnack's theorem states that such curves cannot have more than $\frac{(d-1)(d-2)}{2}+1$ components, but the possible arrangements of their ovals is a far more complicated problem; we refer the reader to the beautiful survey \cite{Wilson} on the subject.

In fact the problem itself of constructing curves with many ovals is a difficult one. Gayet and Welschinger's theorem \cite{GayetWelschinger2011} on the exponential rarefaction of maximal curves confirms such difficulty: if we write the equation of a curve putting some coefficients in front of the standard monomials, 
it is exponentially difficult to get one with many ovals.

This ``approach" of building curves depends on the choice of the basis we want to use for $W_{2,d}$; if instead of the monomial basis (which is orthogonal for the Bombieri scalar product, inducing the Kostlan distribution) we use the spherical harmonics the result improves dramatically. In fact now, because of Theorem \ref{thm:main} for $\alpha=0$ and $n=2$ the expectation of the number of components of the curve we obtain is a fraction of the maximal one. 
The price we pay is that we are no longer allowed to use the simple monomial basis and we need an explicit expression for the polynomials representing spherical harmonics, such an expression being quite complicated.

A completely analogous discussion holds for real hypersurfaces: now the fact that what we get applying Theorem \ref{thm:main} is on average a fraction of the maximal number of components is due to Milnor's bound.

The method employed here (based on the ideas of Nazarov and Sodin) suggests that more generally one could ask for the expected \emph{arrangement} of such ovals. 
In fact essentially the proof relies on the following observation: 
given any point $x$ on the sphere and a small disk around it of radius $\Theta(d^{-1})$ there is a positive probability (independent of $d$) of finding an oval inside such small disk (and such oval loops around $x$). 
Unfortunately this probability, even if constant, is still very little. 
A numerical computation of the expectation of the number of components for $n=2$ was performed by M. Nastasescu in \cite{Nastasescu}: 
as Sarnak writes in his letter ``a random curve is 4\% Harnack".
Thus, the random version of Hilbert's Sixteenth Problem, to study the expected \emph{arrangement} of a random curve, is still open and despite its difficulty it seems to be very fascinating.

One of the main difficulties in this framework is to define a statistic (beyond the number of components) that captures some features of the arrangement. 
We suggest one that we call the \emph{energy} of the curve $C$.

\begin{figure}[ht]
    \begin{center}
    \includegraphics[scale=.3]{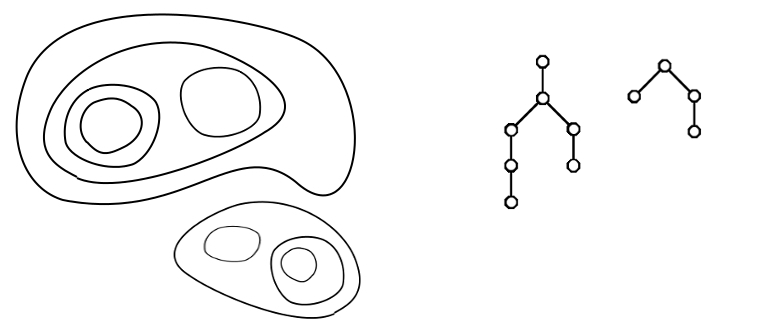}
    \end{center}
    \caption{A hypothetical arrangement of ovals for $C$ and the corresponding forest.  The energy is $h(C) = 2 \cdot 2 \cdot (2\cdot 2 + 2) + 2 \cdot (2 + 2 \cdot 2) = 36$.}
    \label{fig:Tree}
\end{figure}

First we say that one \emph{oval} $C_1$ is contained into another one $C_2$ if it lies in the interior of $C_2$. Using containment we can give a partial ordering on ovals
and construct a forest (a collection of trees) representing their arrangement.
The vertices of the trees represent the ovals, with the empty ovals represented by leaves,
and an oval not contained in any other is represented by a root.
We define an energy that is additive for disjoint unions of trees, and multiplicative when appending trees.
So the energy of a forest is the sum of the energies of its trees.
We define the energy of a tree with only one isolated vertex as $2$,
and of a larger tree inductively as twice the sum of the energies of all subtrees that result from removing its root.
This description determines the energy of any forest (for example, see Figure \ref{fig:Tree}).

By B\'ezout's Theorem the length of each tree cannot be more than $\lfloor \frac{d}{2}\rfloor$ and if such a tree exists then no other trees can appear; 
the corresponding curve is said to be maximally nested and its energy is $h(C)=2^{d/2}$.  On a somehow opposite direction one can consider Harnack curves, 
having $h(C) \sim 5d^2$ (for $d$ large enough). Such curves have the maximal number of components and satisfy some special condition on the arrangement of their ovals; 
we refer the reader to \cite{Mikhalkin} for the exact definition and to \cite{Kenyon} for surprising connections with planar dimers.

Thus, using this notion, one open problem that we can state on the expected arrangement is to describe the asymptotic behavior as $d$ goes to infinity of the energy of a random curve.
The pictures we get of a random curve suggest that the energy has the same order as a Harnack curve,
so we conjecture that: 
$$\lim_{d \rightarrow \infty} \frac{\log \mathbb{E}h(C)}{\log d} = 2$$
(the lower bound follows from Theorem \ref{thm:main}).

\subsection{Plan of the paper}

In Section 2, we specialize to the case $n=1$.
This helps to illustrate the main techniques used in the paper and provides a proof for the value of the parameter $\delta$
used Section 6.2.
Theorem \ref{volume} is proved in Section3 and
Theorem \ref{thm:main} in Section 4.
Section 5 is devoted to the proof of the main technical lemmas for Theorem \ref{thm:main}.
In Section 6 we discuss some examples (more on random curves and random surfaces).
In the Appendix we provide the proof of Milnor's bound,
we review a construction of hyperspherical harmonics,
as well as some of their properties.

 \section*{Acknowledgements}
 The authors are grateful to S. Basu for his constant support and many useful suggestions and A. Lancichinetti for stimulating discussions.
 \section{A toy model proof: the univariate case}

To illustrate the main ideas, we first consider a random polynomial in the special case $n=1$.
As mentioned in the introduction, in this case both the zero dimensional volume and the number of connected components coincide. For simplicity of notation we will denote their expectation by:
$$E_d=\EE\{\textrm{number of zeroes of a real random univariate polynomial}\}.$$
In view of the generalization to higher dimensions, we will give two proofs that $E_d=\Theta(d)$. The first one is based on techniques coming from \emph{integral geometry} and will produce an exact result; this proof will generalize to any number of variables producing again an exact result. The second is based on the \emph{barrier method} introduced by Nazarov and Sodin; it will produce only an approximate result but has the advantage that it generalizes to higher dimensions.

We notice that in the case $n=1$, a real random polynomial is just a trigonometric polynomial with either only even or only odd terms (depending on whether $d$ is even or odd).

\subsection{Integral geometry proof}
The first proof is obtained following the ideas of \cite{EdelmanKostlan95}. 
Recall that if we have functions $\{f_1, \ldots, f_k\}$ defined on the real line we can consider the random function:
$$f=\sum_{i=1}^k \xi_if_i,\quad\textrm{with}\quad \xi_i\sim N(0,1).$$
Then if the coefficients are independently distributed, the expected number of zeros of $f$ on the interval $[a,b]$ is given by the formula:
\begin{equation}\label{expintgeom}\EE\{\textrm{number of zeroes of $f$ on $[a,b]$}\}=\frac{1}{\pi}\int_a^b\| \dot{\gamma}(t) \|dt,\end{equation}
where $\gamma$ is the projection of the curve $c:t\mapsto(f_1(t),\ldots, f_k(t))$ on the unit sphere $S^{k-1}$ (we refer the reader to the beautiful paper \cite{EdelmanKostlan95} for this and related topics). 
We note that (\ref{expintgeom}) is proved using a special case of the \emph{kinematic formula} from integral geometry.
In our case $[a,b]=[0, 2\pi]$ and if $d$ is even,
$$f(t)=\xi_0\frac{1}{\sqrt{2\pi}}+\sum_{k=1}^{d/2} \bigg(a_k \frac{\sin(2kt)}{\sqrt{\pi}}+b_k\frac{ \cos(2kt)}{\sqrt{\pi}}\bigg),\quad \xi_0, a_k, b_k \sim N(0,1).$$
The curve $c$ is given by:
$$c(t)=\frac{1}{\sqrt{\pi}}\bigg(\frac{1}{\sqrt{2}}, \sin(2t), \cos(2t), \ldots, \sin(td), \cos(td)\bigg)$$
and its norm in $\RR^k$ is constant and equals $\sqrt{\frac{d+1}{2 \pi}}.$
In particular we have:
$$\gamma(t) = c(t)\sqrt{\frac{2\pi}{d+1}} \quad \textrm{and}\quad \| \dot{\gamma} \| = \|\dot{c}\| \sqrt{\frac{2\pi}{d+1}}.$$
The norm of $\dot{\gamma}$ is
$$\|\dot{\gamma} \|=\|(0, 2\cos(2t), -2\sin(2t), \ldots, d \cos(td), -d \sin (td))\|\sqrt{\frac{2}{d+1}},$$
which easily simplifies to:
$$\|\dot{\gamma} \|=\sqrt{\frac{ d(d+2)}{3}}.$$
A similar calculation in the case $d$ is odd gives the exact same result.
We are now in the position to use formula (\ref{expintgeom}) and we obtain that the expected number of zeros of $f$ on $[0, 2\pi]$ is $2\sqrt{\frac{ d(d+2)}{3}}$.
Thus, the number $E_d$ of expected zeros on the \emph{real projective line} is:
$$E_d=\sqrt{\frac{ d(d+2)}{3}}.$$
Notice that this is an exact result and gives the asymptotic $E_d=\Theta(d).$
\subsection{Barrier method proof}
In this section we only prove that $E_d\geq c d$ for $d$ large enough; the inequality $E_d\leq d$ follows from the Fundamental Theorem of Algebra.
The proof is centered around the function described in the following Claim whose proof is deferred until the end.

\begin{claim}[Existence of the barrier]
There exist numerical constants $\rho$ and $c_1$ such that, 
for each sufficiently large $d$ and each
$x \in S^1$ 
there is a homogeneous polynomial $b_x $ of degree $d$ and $L^2$-norm one satisfying:
$$b_x(x) \geq c_1 \sqrt{d} \quad \textrm{and}\quad b_x( e^{i \rho / (d+1)} \cdot x) \leq - c_1 \sqrt{d}. $$
\end{claim}

Once a point $x$ in $S^1$ is fixed,
we can decompose $f = \xi_0 b_x + f_x$
where $\xi_0$ is a Gaussian random variable with variance
$\frac{1}{d+1}$, and $f_x$ is
in the orthogonal complement to $b_x$ 
in $W_{1,d}$. We choose a new  Gaussian
random variable $\tilde{\xi}_0$ independent of $\xi_0$ and of $f_x$ 
with variance $\frac{1}{d+1}$,
and set:
$$f_\pm = \pm \tilde{\xi}_0 b_x + f_x .$$ 
By construction the random functions $f_+$ and $f_-$ have the same distribution as $f$. 
We notice that we can write:
$$ f = \xi_0 b_x + \frac{1}{2} \left( f_+ + f_- \right) .$$
Next we apply the following Claim to $f_+$ and $f_-$.
\begin{claim}
There exists a constant $c_2$ such that $\EE |f(x)|\leq c_2$ for every $x$ in $S^1.$

\end{claim}

We also defer the proof of this Claim in favor of first seeing how it is used. Let us recall that Markov's inequality for a positive random variable $|X|$ reads $\PP\{|X|\geq c\}\leq c^{-1}\EE |X|$. If we apply such inequality to $|X|=|f(x)|$, because of the previous Claim we can find a constant $c_0$ such that:
$$\PP \{ |f_{\pm}(x)| \geq c_0 \} \leq 1/5\quad \textrm{and}\quad \PP \{ |f_{\pm}( e^{i r / (d+1)} x)| \geq C_0 \} \leq 1/5,$$
where $\rho$ is provided by the \emph{Existence of the barrier} Claim.

The event that $|f_\pm(x)|$ and $|f_\pm(e^{i \rho / (d+1)} x)|$ are each simultaneously \emph{at most} $c_0$ is
an intersection of four events each of which is complementary to the one of the above.
Thus the probability of their intersection satisfies:
$$\PP \left( \cap \right)  = 1 - \PP \left( \cup \right) \geq 1- 4 \cdot 1/5 = 1/5. $$
We note that these four events are not mutually independent, but this is not needed in the line above.
Now consider the event $\Omega(x)$ that 
$f (x) \geq c_0$ and $f( e^{i \rho / (d+1)} x) \leq - c_0 $. If $\Omega(x)$ occurs then $f$ has a zero\footnote{In this one dimensional case the existence of a zero is essentially the \emph{intermediate-value} property. In higher dimensions this event will be generalized by considering the \emph{boundary} of a disk centered at $x$; in the current situation we only need to use one of the two boundary points of an interval centered at $x$.} on $S^1$ between $x$ and $e^{i \rho / (d+1)} x$.
The event $\Omega(x)$ happens provided that: (i) $\xi_0 \sqrt{d} \geq 2 c_1^{-1} c_0$
and (ii) $f_\pm$ evaluated at at $x$ and $e^{i \rho / (d+1)} x$ are each bounded in absolute value by $c_0$.
Note that (i) and (ii) are independent, and we have just checked above that (ii) happens with probability at least $1/5$.
Moreover, $\PP \{ \xi_0 \sqrt{d} \geq 2 c_1^{-1} c_0 \}\geq c_3,$
for some constant $c_3$, since $\xi_0 \sqrt{d}$ has variance $\frac{d}{d+1}=\Theta(1)$.

It remains to choose $\Theta(d)$ disjoint intervals in $S^1$ each of length $\rho/d$. Each of
them contains a zero of $f$ with probability at least $c_3/5$ and hence: $$E_d \geq c \cdot d ,$$
as desired.

\begin{proof}[First Claim]

We will choose our barrier in such a way that the barrier $b_x$
centered at some point $x \in S^1$ is obtained just by precomposing $b_{(0,1)}$
with a rotation that moves the point $(0,1)$ to $x$,
so we only need to define $b_{(0,1)}$ which we do by normalizing the function: 
$$U_d(e^{i \theta}) = U_d(\cos \theta, \sin \theta) = \frac{\sin( (d+1) \theta )}{\sin \theta}.$$
Notice that $U_d$ is a polynomial of degree $d$, since the denominator divides the numerator.
($U_d$ is a projectivized Chebychev polynomial of the second kind).\\
We have:
$$U_d(0,1) = d ,$$ 
say by de l'Hopital's rule with $\theta \rightarrow 0$, 
and for all $d$ sufficiently large, we also have 
$$U_d \left( e^{i \frac{3 \pi}{2 (d+1)}} \right) \leq - \frac{1}{2} d .$$
Thus it suffices to show that the norm of $U_d$ is on the order of $\sqrt{d}$ as $d \rightarrow \infty$,
for then we will obtain the desired properties for $b_{(1,0)}$ (the normalization of $U_d$).

It is easy to show that $ || U_d ||^2 = 2 \pi d;$ thus the function $$b_{(0,1)}(x) = \frac{1}{\sqrt{2 \pi d}}U_d$$ has the desired properties with 
$c_1 = \frac{1}{2\sqrt{2 \pi}}$ and $\rho = \frac{3\pi}{2(d+1)}$.

\end{proof}

\begin{proof}[Second Claim]

Let us choose a new orthonormal basis $g_0,g_1,..,g_d$ so that all basis elements except $g_0$ vanish at $x$.
To see that it is possible to do this, first notice that the subspace $V^{\perp}$ of homogeneous polynomials of degree $d$
vanishing at $x$ has codimension one inside $W_{1,d}$.
Thus, to get the desired orthonormal basis, 
we may first orthonormalize a basis for $V^{\perp}$, then take its orthogonal complement $V$ which consists of just one function.
Once normalized, this function gives the basis element $g_0$.

Writing our random polynomial in terms of this basis,
we have
$$|f(x)| = \big| \xi_0 g_0(x) + \sum_k \xi_k g_k(x)\big| = |\xi_0| \cdot |g_0(x)|.$$
Using again the standard basis $f_k$, the fact that $g_0$ is a normalized trigonometric polynomial of degree $d$ means that 
$$|g_0(x)| = \left| \sum_k A_k f_k(x) \right| \leq \sum_k |A_k| \leq \sqrt{\sum_k |A_k|^2}\sqrt{d+1} = \sqrt{d+1},$$
where we have also used the triangle inequality and the discrete Cauchy-Schwarz inequality.
Combining this with the above,
$$|f(x)| \leq | \sqrt{d+1}\cdot \xi_0|,$$
which implies that $\EE |f(x)|$ is bounded by a constant (the Gaussian random variable $\sqrt{d+1}\cdot \xi_0$ has variance of constant order,
and therefore so has the expectation of its modulus).

\end{proof}

 \section{Proof of Theorem \ref{volume}}
We generalize the proof given in \cite{Berard} for the \emph{harmonic} case to our $\alpha$-window of eigenspaces. 
Here the main idea is to use a basis of the eigenfunctions of the Laplacian to define a kind of \emph{Veronese} embedding, 
transforming the geometry of zeros of such functions on $S^n$ to the geometry of the intersection of the image of $S^n$ with hyperplanes. 
What is nice about using harmonic functions is that such a Veronese embedding turns out to be a dilation, i.e. dilates the Riemannian metric by a constant factor.  
In the language of integral geometry such embeddings are called \emph{moment maps} (generalizing the idea of the \emph{moment curve} $\gamma$ that we used in the univariate case
in Section 2.1).

 Let us denote by  $d(n,l)$ the dimension of $H_{n,l}$ as above (thus in the case $n=2$ we have $d(l,2)=2l+1$ and in the general case $d(n,l)=\Theta(l^{n-1})$). 
 For each $(n,l)$ let $\{Y_i\}_{i=1}^{d(n,l)}$ be an orthonormal basis for $H_{n,l}$ and consider the map:
 $$\Lambda_{n,l}:S^n\to \RR^{d(n,l)}$$
 whose components are $(Y_1, \ldots, Y_{d(n,l)}).$ Such a map has the following properties: 
 (i) it is an immersion and its image is a submanifold of the sphere in $\RR^{d(n,l)}$ of radius $R=(d(n,l)/\textrm{Vol}(S^n))^{1/2}$; 
 (ii) it is a dilation of ratio $\lambda R^2/n$ where $\lambda=l(l+n-1)$ it is the eigenvalue for the eigenspace $H_{n,l}$ of $\Delta_{S^n}$. 
 In other words the length of the image of a unit vector in $S^n$ under the differential of $\Lambda_{n,l}$ is $\sqrt{\lambda/n}R$. 
 Property (i) follows from the addition theorem for spherical harmonics:
 $$Z_l(x, y)=\sum_{i=1}^{d(n,l)}Y_i(x)Y_i(y),\quad x,y\in S^n.$$
 Here $Z_l$ denotes the (reproducing kernel) zonal harmonic of degree $l$,
 for which  $Z_l(x,x) = \frac{d(n,l)}{\textrm{Vol}(S^n)}$.
 Thus, evaluating the addition formula at $y=x$ we get:
 $$\frac{d(n,l)}{\textrm{Vol}(S^n)}=\|\Lambda_{n,l}(x)\|^2,$$
 which shows the image is contained in the sphere. Property (ii) is proved in \cite{Besse} and goes under the name of \emph{nice imbeddings property}.
 
 For a given $\alpha \in [0, 1)$ we consider thus the map:
 $$\Lambda:S^n\to \RR^{D_\alpha}=\bigoplus_{l\in [\alpha d, d]}\RR^{d(n,l)}$$
 whose components are the $\Lambda_{n,l}$ normalized by the factor $(D_\alpha/\textrm{Vol}(S^n))^{1/2}$, i.e. 
 $$\Lambda=\bigg(\frac{\textrm{Vol}(S^n)}{D_\alpha}\bigg)^{1/2}(\Lambda_{\alpha d}, \ldots, \Lambda_l, \ldots, \Lambda_d), \quad d-l\textrm{ even}.$$
 With this choice of the normalization the addition theorem for spherical harmonics implies that the image of $\Lambda$ is contained into the unit sphere in $\RR^{D_\alpha}$. We claim now that $\Lambda$ is a dilation of ratio $\delta$, as defined in the Introduction. To prove this let us consider a unit vector $v\in TS^n$ and let us compute the square of the norm of its image under the differential of $\Lambda$:
 $$\|d\Lambda (v)\|^2=\sum_{l\in [\alpha d, d]} \frac{\textrm{Vol}(S^n)}{D_\alpha}\|d\Lambda_{n,l} (v)\|^2=\sum_{l\in [\alpha d, d]}\frac{\lambda (n,l) d(n,l)}{n D_\alpha}=\delta,$$
 where in the second equality we have used the fact that $\Lambda_{n,l}$ is a dilation of ratio $\lambda(n,l)R^2/n$. We apply now the integral geometry formula (see \cite{EdelmanKostlan95}):
 \begin{equation}\label{intgeom}\frac{1}{\textrm{Vol}(S^{D_\alpha -1})}\int_{S^{D_\alpha -1}}\textrm{Vol}(\{a\}^\perp \cap \Lambda(S^n))da=\textrm{Vol}(\Lambda(S^n))\frac{\textrm{Vol}(S^{n-1})}{\textrm{Vol}(S^n)}.\end{equation}
 Property (ii) above implies that the map $\Lambda$ is a covering space (in fact either a double cover, or a diffeomorphism); let us say it is $k$-sheeted. We claim now that the left hand side of (\ref{intgeom}) equals exactly $k\delta^{(n-1)/2}\EE \textrm{Vol}(Z_{S^n}(f))$. In fact if $a=(a_{1}, \ldots ,a_{D_\alpha})$, then $\{a\}^\perp \cap \Lambda (S^n)$ equals the image under $\Lambda$ of the nodal line of the function $\sum_{i=1}^{D_\alpha} a_i Y_i$, where $\{Y_i\}_{i=1}^{D_\alpha}$ is the basis for $H_{[\alpha d, d]}$ obtained by collecting together the elements from the various orthonormal basis of $H_{n,l}$. Since $\Lambda$ is a covering, the volume of such image (which is $n-1$ dimensional) equals $k \delta^{(n-1)/2}$ times the original volume; moreover since our Gaussian distribution is uniform on the unit sphere in $H_{[\alpha d, d]}$, then the above integral computes expectation and the conclusion follows.
 
 As for the right hand side of (\ref{intgeom}) again the fact that $\Lambda$ is a $k$-sheeted covering and a dilation of ratio $\delta$, implies $\textrm{Vol}(\Lambda(S^n))=k\delta^{n/2}\textrm{Vol}(S^n)$. Plugging all this into (\ref{intgeom}) we finally get:
 $$\EE \textrm{Vol}(Z_{S^n}(f))=\delta^{1/2}\textrm{Vol}(S^{n-1}).$$

\section{Proof of Theorem \ref{thm:main}}
Our proof is deeply influenced by the ideas and the exposition of \cite{NazarovSodin2009}. In fact for $\alpha=1$ our claim is the direct generalization of the \emph{lower bound estimate} in \cite{NazarovSodin2009}; 
this case is a simple extension of the proof in \cite{NazarovSodin2009} and we devote a remark to it at the end of the proof.

Notice that the upper bound for our expectation is provided by Milnor's bound (\ref{Milnor}), discussed in the Appendix, and we only have to prove the lower bound.

Thus we fix\footnote{In particular every constant in this proof is allowed to depend on $\alpha$, but not on $d.$} for the rest of the proof
 an $\alpha\in [0,1).$ To start with, for $x$ a point on the sphere $S^n$ we consider the following event:
\begin{equation}\label{omega}\Omega(x, r)=\{\f(x)>0\quad \textrm{and} \quad \f |_{\partial D(x, r)}<0\}.\end{equation}
If $\Omega(x, r)$ occurs, then $\{\f=0\}$ has a component inside $D(x, r)$ (and this component \emph{loops} around $x$). If we prove that for $r=\Theta(d^{-1})$ there exists a constant $a_1>0$ independent of $x$ and $d$ such that
$$\PP\{\Omega(x, r)\}\geq a_1,$$
then we can cover the sphere $S^n$ with $\Theta(d^n)$ disjoint  such disks and get the statement (each one of these disks contributing at least $a_1$ to the expectation in the statement).

Specifically, the radius $r=r(d,n)$ we are going to use from now on is given by:
$$r=\frac{2y_n}{2d+n-1},$$
where $y_2$ is the first point of minimum of the Bessel function $J_0$ and $y_n$ for $n>2$ is the first point \emph{to the right of zero} where the Bessel function $J_{\frac{n-2}{2}}$ reaches a minimum (in this case $J_{\frac{n-2}{2}}$ already has a minimum at zero, see \cite{Szego}).

The following Lemma provides the existence of the so called \emph{barrier} function; we defer its proof until the next section.
\begin{lemma}[Existence of a barrier] \label{lemma1}
For every $\alpha \in [0,1)$ there exists a constant $c_1=c_1(\alpha)>0$ such that for every $d>0$ and every $x$ in $S^n$ there exists a homogeneous polynomial $B_x$ of degree $d$ and $L^2$-norm one satisfying:
$$B_x(x)\geq c_1 d^{n/2} \quad \textrm{and}\quad B_x|_{\partial D(x, r)}\leq -c_1 d^{n/2}.$$
\end{lemma}

\begin{figure}[ht]
    \begin{center}
    \includegraphics[scale=.4]{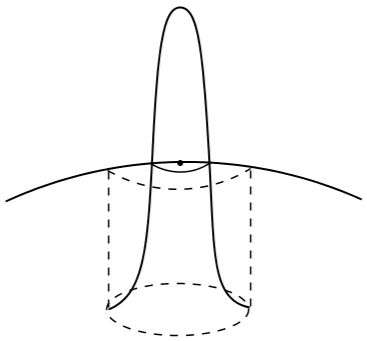}
    \end{center}
    \caption{The behavior of the barrier function $B_x$ at the point $x$.}
    \label{fig:barrier}
\end{figure}

Using $\B$ we can decompose $\HH=\textrm{span}\{\B\}\oplus\textrm{span}\{\B\}^\perp$, 
thus getting the following decomposition for $\f$:
$$\f=\xi_0 \B+\f^\perp\quad\textrm{with}\quad \xi_0\sim N(0, D_\alpha^{-1}).$$
We let now $\tilde \xi_0$ be a random variable distributed as $\xi_0$ but independent of it and we define ${\f}_\pm=\pm \tilde \xi_0 \B+\f^\perp$. 
Notice that both $\f$ and ${\f}_\pm$ have the same distribution. The introduction of these new random polynomials allows us to write:
\begin{equation}\label{splitf}\f=\xi_0 \B+\frac{1}{2}({\f}_{+}+{\f}_{-})\end{equation}
and to split our problem into the study of the behavior of $\B$ and ${\f}_\pm$ separately. 
In fact the event $\Omega(x, r)$ happens provided that for some constant $a_2>0$ the two following events both happen:
\begin{itemize}
\item[1)]$E(x,r)=\{\xi_0\B(x)\geq 2a_2 \textrm{ and } \xi_0B_x|_{\partial D(x, r)}\leq -2a_2\}$;
\item[2)]$G(x,  r)=\{\|{\f}_\pm(x)\|\leq a_2$ and $\|{\f}_{\pm}|_{\partial D(x, r)}\|\leq a_2\}.$
\end{itemize}
To check that $E(x,r)\cap G(x,r)$ implies $\Omega(x, r)$ we simply substitute the inequalities defining $E(x,r)$ and $G(x, r)$ in (\ref{splitf}) evaluating respectively at $x$ and at $\partial D(x, r)$.

Now by definition the two events $E(x,r)$ and $G(x, r)$ are independent and thus:
$$\PP\{E(x,r)\cap G(x, r)\}=\PP\{E(x,r)\}\PP\{G(x, r)\}.$$
Because of this it suffices to show that the probability of each one of them is bounded from below by a positive constant that does not depend on $d$ and this will provide a bound from below for the probability of $\Omega(x, r)$ which is independent of $d$.

The fact that the probability of $E(x,r)$ is bounded from below independently on $d$ immediately follows from Lemma \ref{lemma1}: since $D_\alpha=\Theta(d^n)$ and $B_x(x)=\Theta(d^{n/2})$ the random variable $\xi_0\B(x)$ is Gaussian with mean zero and variance $\Theta(1)$; similarly for the boundary. Thus for $a_3$ small enough the probability of it being bigger than $a_3$ is uniformly bounded from below.

It remains to study $G(x, r)$; to this extent notice that this event is given by the intersection of the four events $G_{1,2}=\{\|{\f}_\pm(x)\|\leq a_2\}$ and  $G_{3,4}=\{\|{\f}_\pm|_{\partial D(x, r)}\|\leq a_2\}$. For these we do not need independence as:
$$\PP\big\{\bigcap_{i} G_i\big\}=1-\PP\big\{\big(\bigcap_{i} G_i\big)^c\big\}=1-\PP\big\{\bigcup_{i} G_i^c\big\}\geq 1-\sum_i \PP\{G_i^c\}$$
and it is enough to prove that the probability of the complement of each one of them is small.

To bound the probability of $G_{1}^c$ and $G_2^c$ we argue as follows. First notice that since $f_\pm$ are distributed as $f$ it is enough to prove the corresponding bound for a random $f$. 
We consider in  $H_{[\alpha d, d]}$ the hyperplane $V^\perp$ of those functions vanishing at $x$ together with its orthogonal complement $V$. 
Such orthogonal complement \emph{must} be spanned by a norm-one function $F$ of the form:
\begin{equation}\label{bestbarrier}F=F_0+\sum_{l\in [\alpha d, d]}y_lY_l\quad \textrm{with}\quad \sum_{l\in [\alpha d, d]}y_l^2\leq 1,\end{equation}
where $Y_l\in H_{n,l}$ is the normalized \emph{zonal harmonic} centered at $x$ (notice in particular that this normalization produced the r.h.s. inequality in the above line).  
In fact in each space $H_{n,l}$ the orthogonal complement of the zonal $Y_l$ centered at $x$ vanishes at this point, and we immediately get (\ref{bestbarrier}).
Such zonal function is obtained by composing the function in equation (\ref{zonal}) with a rotation that moves the north-pole of $S^n$ to $x$; 
moreover by setting $\theta=0$ in equation (\ref{zonal}) we see that it has order $\Theta(l^{(n-1)/2})$ at $x$. 
Using the above orthogonal decomposition $H_{[\alpha d, d]}= V\oplus V^\perp$, 
we can write the random $f$ as $\xi_0 F+F^\perp$ with $F^\perp$ a random function in $V^\perp.$ Thus we get:
\begin{equation}\label{bestbarr}\EE |f(x)|\leq \EE \bigg|\xi_0\sum_{l\in [\alpha d, d]}|y_lY_l(x)|\bigg|\leq \EE\bigg|\xi_0\bigg(\sum_{l\in [\alpha d, d]}Y_l(x)^2\bigg)^{1/2}\bigg|\leq a_3.\end{equation}
We have used (\ref{bestbarrier}) and the discrete Cauchy-Schwarz inequality in the second inequality; for the last inequality we have used the fact that $Y_l(x)^2=\Theta(l^{n-1})$, which implies
$$s(d)=\bigg(\sum_{l\in [\alpha d, d]}Y_l(x)^2\bigg)^{1/2}=\Theta(d^{n/2}).$$ Since the variance of $\xi$ has order $\Theta(d^n)$, when we scale it by $s(d)$ we make its variance of constant order, which gives (\ref{bestbarr}). The statement on the probability of $\{|f(x)|\geq a_2\}$ being smaller immediately follows from Markov's inequality.

We explain now the argument for the remaining two cases. Again since ${\f}_\pm$ is distributed as $\f$, the probability of each $G_{3,4}^c$ is exactly the probability of the following event:
$$\{\|\f|_{\partial D(x, r)}\|\geq a_2\quad \textrm{for a random $\f \in \HH$}\}.$$
Hence we consider the positive random variable: $$|X|=\max_{\partial D(x, r)}|\f|.$$
Bounding the expectation of $|X|$ is more delicate and we devote to it a separate lemma, whose proof is postponed to the next section.
\begin{lemma}\label{lemma2} 
For every $\alpha\in [0,1]$ there exists a constant $c_2=c_2(\alpha)$ such that 
$$\mathbb{E} \max_{\partial D(x, r)}|f|\leq c_2.$$
\end{lemma}
Applying Markov's inequality to $|X|$ combined with the previous Lemma gives a bound for the probability of $G_{3,4}^c,$ concluding the proof.

\begin{remark}
This proof does not apply when $\alpha=1$, but in that case
one can use a rather direct generalization to $n+1$ variables of the proof given in \cite{NazarovSodin2009}.
Namely, their Claim 3.1 providing a bound on the expected max over a small disk 
is still true in more variables and can be proved in the same way,
and the barrier function can again be chosen to be the normalized zonal harmonic composed with a rotation.
\end{remark}

\section{Proofs of Lemmas}
\subsection{Proof of Lemma \ref{lemma1}}\label{sec:lemma1}
We work out the details for $n=2$ first.
We will construct $B_x$ for $x$ the north pole and define it for the other points by considering the composition with rotations (such rotations preserve all the properties of the statement).

Set $y=y_2$ (the first minimum point of the Bessel function $J_0$) and let $Y_l^0$ be the normalized \emph{zonal} harmonic of degree $l$ centered at the North Nole. Using the \emph{Legendre polynomial} $P_l$ of degree $l$ we can write:
$$Y_l^0(\theta, \phi)=P_l(\cos \theta)\sqrt{2l+1},$$where $(\theta, \phi)$ are the standard polar coordinates on $S^2$. We recall from \cite{Szego} the following asymptotic for $P_l$ (Hilb's asymptotic):
$$P_l(\cos \theta)=\bigg(\frac{\theta}{\sin \theta}\bigg)^{\frac{1}{2}}J_0\bigg(\frac{2l+1}{2}\theta\bigg)+R_l(\theta),$$
where $J_0$ is the Bessel function of the first kind and order zero. The error $R_l(\theta)$ is always less than $\theta^{1/2}O(l^{-\frac{3}{2}}).$
The function $J_0$ satisfies the following properties: $J_0(0)=1,$ it is decreasing in $[0,y]$ where $y$ is a local minimum with $J_0(y)<0$. 
We call $z\in (0,y)$ the first zero of $J_0$. Let $\mu$ be the maximum of the two values $\alpha$ and $z/y$; we define $B_x$ as:
$$B_x=\frac{\sum_{l\in [\mu d, d]}Y_l^0}{\sqrt{\textrm{Card}[\mu d, d]}}.$$
Since each $Y_l^0$ has norm one and since different zonals are pairwise orthogonal, we immediately see that the $L^2$-norm of $B_x$ is \emph{one}, as required from the statement. 
Let us now estimate the value of $B_x$ at $x$. Using the integral estimate for the numerator, which provides a bound of the form $\Theta(d^{3/2})$, 
and the fact that $\sqrt{\textrm{Card}[\mu d, d]}=\Theta(d^{1/2})$, for some constant $a_1>0$ independent of $d$ we have:
$$B_x(x)=\frac{\sum_{l\in [\mu d, d]}\sqrt{2l+1}}{\sqrt{\textrm{Card}[\mu d, d]}}\geq a_1 d.$$ 
Recall that the number $r$ was defined such that $\frac{2d+1}{2}r=y,$ i.e. we set:
$$r=\frac{2y}{2d+1}.$$
We estimate now the value of $B_x$ on the boundary of a disk of radius $r$. Notice that it suffices to show the bound for the numerator in the definition of $B_x$ has order
$\Theta(d^{3/2})$, since we already know the denominator has order $\Theta(d^{1/2})$. Because of Hilb's asymptotic we have:
$$\Theta (d^{1/2})B_x(r)=\underbrace{\sum_{l\in [\mu d, d]}\sqrt{\frac{r}{\sin r}}J_0\bigg(\frac{2l+1}{2}r\bigg)\sqrt{2l+1}}_{L(r)}+\underbrace{ \sum_{l\in [\mu d, d]}R_l(r)\sqrt{2l+1}}_{R(r)}.$$
We estimate the error $R(r)$ first. 
Notice that the summands $R_l(r)\sqrt{2l+1}$ are all bounded by $r^{1/2} O(l^{-1})$
so that: 
$$R(r) \leq r^{1/2} O \left( \sum_{l=1}^d l^{-1} \right) = r^{1/2} O(\log d) = O \left( \frac{\log d}{\sqrt{d}} \right).$$

Next we bound $L(r)$ as follows. First notice that since $\mu=\max\{\alpha, z/y\}$, then for every $l$ in the sum $J_0(\frac{2l+1}{2}r)<0$. 
We fix now a $\gamma\in (\mu, 1)$ and notice that since each term is negative and $J_0$ is decreasing in $(\gamma y, y)$, then:
$$L(r)\leq \sum_{l\in [\gamma d, d]}\bigg(\frac{r}{\sin r}\bigg)^{1/2}J_0(\gamma y)\sqrt{2l+1}.$$
Moreover using the expansion $(r/\sin r)^{1/2}=1+\Theta(d^{-1})$ we get:
$$L(r) \leq J_0(\gamma y)(1+\Theta(d^{-1}))O(d^{3/2})\leq -c_2d^{3/2},$$
where we have used the integral estimate for the sum of the square roots. 
Recalling that $L(r)=B_x(r)\Theta(d^{1/2})$, this concludes the proof for $n=2$.

For the case $n>2$ we argue as follows. First recall that we set $y_n$ to be the first point after zero where the Bessel function $J_{\frac{n-2}{2}}$ reaches a minimum. 
The function $J_{\frac{n-2}{2}}$ has the following properties: it vanishes at zero and we denote by $x_n$ its first point of maximum (it is positive at this point); it is decreasing on  the interval $[x_n, y_n]$ and negative at $y_n$. We let $z_n$ be its only zero in $[x_n, y_n]$ and $\mu_n$ be $\max \{\alpha, z_n/y_n\}.$ If $Y_l\in H_{n,l}$ is the normalized zonal harmonic centered at $x$, we define the barrier as:
$$B_x=\frac{\sum_{l\in [\mu_n d, d]}Y_l}{\sqrt{\textrm{Card}[\mu_n d, d]}}.$$
Since the zonals are pairwise orthogonal, 
the $L^2$-norm of $B_x$ is one. 
As already discussed in the proof of Theorem \ref{thm:main}, the order at $x$ of each zonal $Y_l$ is $\Theta(l^{(n-1)/2})$ as shown in (\ref{zonalorder}), and thus the order of their sum at $x$ is $\Theta(d^{(n+1)/2})$. Moreover since $\sqrt{\textrm{Card}[\mu_n d, d]}=\Theta(d^{1/2})$, then we immediately get that for some positive constant $c_1$:
$$B_x(x)\geq c_1 d^{n/2}.$$
The proof of the bound for the value of $B_x$ on $\partial D(x, r)$ involves the same steps as in the case $n=2$, 
this time using a generalization of Hilb's asymptotic involving the Bessel function $\displaystyle J_{\frac{n-2}{2}}$. 
We need an explicit expression for $Y_l$ first:
$$Y_l(\theta)=c(n,l) P_{l}^{(\frac{n-2}{2},\frac{n-2}{2})}(\cos \theta),$$
where $P_{l}^{(\frac{n-2}{2},\frac{n-2}{2})}$ is a Jacobi polynomial
and $c(n,l)$ has order $\Theta(l^{1/2})$, see (\ref{c(n,l)}) in the Appendix. 

Using the asymptotic provided in (\cite{Szego}, Theorem 8.21.12) we have:
$$P_{l}^{(\frac{n-2}{2},\frac{n-2}{2})}(\cos \theta) = \bigg(\sin \frac{\theta}{2}\cdot\cos \frac{\theta}{2}\bigg)^{\frac{2-n}{2}}\bigg\{ h(n,l)\bigg(\frac{\theta}{\sin \theta}\bigg)^{1/2}J_{\frac{n-2}{2}}(N \theta)+R_l(\theta)\bigg\},$$
where $N=\frac{2l + n-1}{2},$ $h(n,l)=\Theta(1)$ is given by (\ref{h(n,l)}), and the error term $R_l(\theta)$ is always less than $\theta^{1/2}O(l^{-\frac{3}{2}}).$

Since by definition $r=\frac{2y_n}{2d+n-1}$ we have $(\sin \frac{r}{2}\cdot\cos \frac{r}{2})^{\frac{2-n}{2}}=\Theta(d^{\frac{n-2}{2}}),$ and plugging all this into the definition of $B_x$, 
we have that $B_x(r)$ equals:
$$\Theta(d^{-\frac{1}{2}})\Theta(d^{\frac{n-2}{2}})\sum_{l\in [\mu_n d, d]}\Theta(l^{1/2})\bigg\{\bigg(\frac{r}{\sin r}\bigg)^{1/2}J_{\frac{n-2}{2}}(r(2l+n-1)/2)+R_l(r)\bigg\}.$$
For the error term we get a bound of the form $O(d^{\frac{n-2}{2}}).$
As for the principal term, we fix as above a $\gamma_n\in (\mu_n, 1)$ and notice that since each summand is negative and $J_{\frac{n-2}{2}}$ is decreasing in $(\gamma_n y_n, y_n)$, then:
$$B_x(r)\leq \Theta(d^{\frac{n-3}{2}})\sum_{l\in [\gamma_nd, d]}\Theta(l^{1/2})J_{\frac{n-2}{2}}(\gamma_n y_n)\leq -c_1d^{n/2}$$
where in the last inequality we have used the integral estimate again.

\subsection{Proof of Lemma \ref{lemma2}}
\label{sec:lemma2}

Since the distribution is invariant by rotation, 
it is enough to prove the statement in the case the point $x$ equals the north-pole $(0,0,..,0,1)$. 
We need an explicit expression for a real orthonormal basis of $H_{[\alpha d,d]}$.
To this end, we use the inductive construction of an orthonormal basis for $H_{n,l}$, the space of hyperspherical harmonics in $n+1$ variables of degree $l$, 
given in terms of orthonormal bases $\{Y_j(\phi)\}_{j \in I_m}$ for each space $H_{n-1,m}$ with $m=0,1,..,l$, where $I_m$ is an index set of size 
$|I_m| = \dim H_{n-1,m} = d(n-1,m)$.
Namely, the basis elements as reviewed in the Appendix are:
$$Y_{l}^m(\theta,\phi) = N_l^m (\sin\theta)^m P^{\left( \frac{n-1}{2}+m \right)}_{l-m}(\cos\theta) Y_j(\phi).$$
In the previous formula $P^{(\lambda)}_{k}$ is the Gegenbauer polynomial \cite[Formula (4.7.1)]{Szego}, $\theta$ is the angle from the north-pole, 
the multi-angle $\phi\in S^{n-1}$  gives the remaining coordinates, and the normalization constant $N_l^m$ is stated explicitly in (\ref{N}) 
(see also \cite[Chapter 3]{Avery} for an exposition of this inductive construction).

Using this basis we can write the random $f$ in $H_{[\alpha d, d]}$ as:
\begin{align*}
f &= \sum_{l\in [\alpha d,d]} \sum_{m=0}^{l} \sum_{j \in I_m}\xi_l^j Y_l^m(\theta,\phi) \\
  &= \sum_{l\in [\alpha d,d]} \sum_{m=0}^{l} \sum_{j \in I_m}\xi_l^j N_l^m (\sin\theta)^m P^{\left( \frac{n-1}{2}+m \right)}_{l-m}(\cos \theta) Y_j(\phi),
\end{align*}
with each $\xi_l^j$ distributed as $N(0,D_{\alpha} ^{-1})$, with $D_\alpha^{-1}=\Theta(d^{-n})$. 
In particular, the restriction of $f$ to $\partial D(x, r)$ is given by:
$$ f(r,\phi) =  \sum_{m=0}^{d} \sum_{j \in I_m} \bigg\{ \sum_{l\in L_m} \xi_l^j N_l^m P^{\left( \frac{n-1}{2}+m \right)}_{l-m}(\cos r)\bigg\} (\sin r)^m Y_j(\phi), $$ 
where the index set $L_m:= \{ l\in[\alpha d,d]\, |\, l \geq m \}$ arises from having reversed the order of summation.

 Notice that the term in the large brackets is a Gaussian random variable itself: it is a linear combination of independent random variables, 
 each one of mean zero and variance $D_\alpha^{-1}$.
 We denote such a random variable by $\hat{\xi}^j$; in other words we define:
 \begin{equation}\label{eq:newxi}
  \hat{\xi}^j =\sum_{l\in L_m} \xi_l^j N_l^m P^{\left( \frac{n-1}{2}+m \right)}_{l-m}(\cos r),\quad j\in I_m.
 \end{equation}
By the addition formula for independent Gaussians we get that $\hat{\xi}^j$ is distributed as a Gaussian with mean zero and variance $\sigma(m,d)^2$ where:
\begin{equation}\label{eq:sigma}
\sigma(m,d)^2 \leq \sum_{l=m}^d \big\{N_l^m P^{\left( \frac{n-1}{2}+m \right)}_{l-m}(\cos r)\big\}^2 .
 \end{equation}
\begin{remark} Most importantly, for each $m$ the Gaussians $\hat{\xi}^j$ are \emph{identically distributed} across all $j$ within the index set $L_m$.
This follows from the fact that in the sum (\ref{eq:newxi}) the coefficient of $\xi_l^j$  depends on $m$ but not $j$.
\end{remark}
Using this new notation we can rewrite:
$$ f(r,\phi) =  \sum_{m=0}^{d} \sum_{j \in I_m} \hat{\xi}^j (\sin r)^m Y_j(\phi) = \sum_{m=0}^{d} (\sin r)^m \sum_{j \in I_m} \hat{\xi}^j Y_j(\phi). $$ 
Using the triangle inequality, we bound the expectation of $|X|=\max_{\partial D(x,r)} |f|$ as:
$$\mathbb{E}|X|\leq \sum_{m=0}^{d} (\sin r)^m \EE \max_{\phi \in S^{n-1}} \left| \sum_{j \in I_m} \hat{\xi}^j Y_j(\phi) \right| .$$
By the previous remark the function:
$$ F(\phi)=\sum_{j \in I_m} \hat{\xi}^j Y_j(\phi)$$
is a random spherical harmonic on $S^{n-1}$ with independent identically distributed Gaussian coefficients. 
This allows us to apply the following Claim, whose proof is provided in the Appendix.

\begin{claim}[A basic estimate for random spherical harmonics:] 
$$\EE \max_{\phi \in S^{n-1}} \left|F(\phi) \right| \leq \sigma(m,d) \cdot d(n-1,m),$$ 
where recall that $d(n-1,m)$
is the dimension of the space of spherical harmonics of degree $m$ in $n$ variables.
\end{claim}

Returning to  $\EE |X|$, we have:
\begin{equation}\label{max}\mathbb{E}|X|\leq \sum_{m=0}^{d} (\sin r)^m \sigma(m,d) \cdot d(n-1,m).\end{equation}
It remains to estimate $\sigma(m,d)$.
 We note that the Gegenbauer polynomials satisfy \cite[Formula (7.33.1)]{Szego}: 
 $$\max_{x\in [-1,1]} \bigg|P^{\left( \frac{n-1}{2}+m \right)}_{l-m}(x)\bigg|=  {l + m + n - 2 \choose l-m}.$$
Using this to estimate (\ref{eq:sigma}) we get:
$$\sigma(m,d)^2 \leq \frac{1}{D_\alpha}\sum_{l=m}^d (N_l^m)^2 {l + m + n - 2 \choose l-m} ^2.$$
Inspecting the constant $(N_l^m)^2$ (see formula (\ref{N})), we notice that we can bound part of it by a constant:
$$ \frac{ \Gamma( \frac{n+2m+1}{2} ) }{ \sqrt{\pi} \Gamma( \frac{n+2m}{2} ) (n+2m-1) }\leq a_0,$$
(in fact this number goes to zero as $m$ goes to infinity); in particular
we bound each term in the above sum by a constant times
$$\frac{ (n+2m-2)! (2l+n-1)(l-m)! }{(l+m+n-2)! }  {l + m + n - 2 \choose l-m} ^2,$$
which simplifies to
$$(2l+n-1) {l + m + n - 2 \choose l-m}.$$
This finally gives:
 \begin{align*}\sigma(m,d)^2& \leq  \frac{a_1}{D_\alpha} \sum_{l=m}^d (2l+n-1) {l + m + n - 2 \choose l-m}\\
 &\leq \frac{2 a_1}{D_\alpha \cdot (2m + n - 2)!}\sum_{l=m}^d(l+m+n-2)^{2m+n-1} \\
 &\leq \frac{2 a_1}{D_\alpha \cdot (2m + n - 2)!} (d+m+n-2)^{2m+n} \\
 &\leq a_2^2 \frac{(2d)^{2m}}{(m!)^2}
 \end{align*}
 where $a_2$ is a constant that does not depend on $d$.
 
In particular $\sigma(m,d) \leq a_2\frac{(2d)^m}{m!}$.
 Recalling that $r=\frac{2y}{2d+n-1}\leq \frac{y}{d}$ (we have set $y_n=y$ for simplicity of notation) 
 and using $(\sin r)^m \leq r^m \leq \left( \frac{y}{d} \right)^m$, we have in (\ref{max}):
 $$\EE |X|\leq \sum_{m=0}^d \left( \frac{y}{d} \right)^m \cdot \sigma(m,d) \cdot d(n-1,m) \leq a_2\sum_{m=0}^d \frac{(2y)^m}{m!} \cdot d(n-1,m).$$
 We estimate $d(n-1,m) \leq a_3 m^{n-2} \leq {a_4}^m$, so that using the definition of the exponential function we obtain:
$$\EE |X| \leq a_5 \sum_{m=0}^d \frac{(2y a_4)^m}{m!}\leq a_5e^{2y a_4} \leq c_2 .$$
 This concludes the proof.

\section{Examples}
\subsection{More on the arrangements of random curves}
Besides the number of ovals of a curve in the projective plane, there are other interesting invariants associated to it that are worth studying (e.g. the \emph{energy} we defined in the Introduction).

One classical invariant is the number of \emph{empty} ovals; let us denote it by $\nu_0$. If the curve $C$ has even degree $d=2k$, then Arnold's inequalities (see \cite{Arnold} section 6) imply:
$$\nu_0(C)\geq b_0(C) -(k-1)(k-2).$$
In particular if a curve has the maximal number of ovals, then $\nu_0(C)\geq k^2.$ The numerical results obtained by M. Nastasescu  in \cite{Nastasescu} have shown that the constant in the bound for $\mathbb{E}b_0$ is very small:
$$\mathbb{E}b_0(C)\leq \frac{1}{25} d^2,$$
and combining this with Arnold's bound would only give $\mathbb{E}\nu_0(C)\geq 0.$ 
On the other hand, as stated in the Introduction, 
the barrier method can be used to give information also on the expected arrangement of the components of the random hypersurface. 
In fact the event $\Omega(x,r)$ as defined above (\ref{omega}) implies the existence of an \emph{empty} oval inside the disk $D(x,r)$.
Thus, placing on the sphere $\Theta(d^2)$ disjoint such disks gives the following:

$$\textrm{Average number of empty ovals:}\quad\mathbb{E}\nu_0(C)=\Theta(d^2).$$

\begin{remark}
We notice that every $f\in W_{n,2k}$ defines a function on the real projective space $\RP^n$ 
and for a generic $f$ the gradient of this function does not vanish on its zero locus. 
In particular if $C$ is a component of $Z_{\RP^n}(f)$, then $\RP^n\backslash C$ has two components; if $q:S^n\to \RP^n$ is the quotient map, we define the \emph{interior} of $C$ to be the component of $\RP^n\backslash C$ whose preimage under $q$ is not connected. 
As we did for curves, we say that a component $C_1$ of $Z_{\RP^n}(f)$ contains another $C_2$ if $C_2$ is contained in the interior of $C_1$; an empty component is one whose interior is connected. 

Using this notation the above corollary can be stated for hypersurfaces as well and provides the existence of $\Theta(d^n)$ empty components.

Similarly we can define the \emph{energy} of an hypersurface as we did for curves: additive on disjoint unions and multiplicative on nestings. 
Instead of always taking the value $2$,
the ``seed'' we take for the energy of an empty component is the sum of its Betti numbers 
(in the case of an empty oval, this equals $2$).\end{remark}

As for higher moments of $b_0$, Nazarov ad Sodin have proved that for $\alpha=1$ (random spherical harmonics) it exponentially concentrates about its mean and Sarnak has informed the authors that the same result holds for $\alpha=0$ (homogeneous polynomials). Thus the \emph{typical} curve has $\Theta(d^2)$ components. The authors expect the same results to hold for any $n$ as well.

On a different direction Gayet and Welschinger have studied random curves in the Kostlan ensemble (the \emph{complex} Fubini-Study one). In \cite{GayetWelschinger2011} they proved that maximal curves become exponentially rare in their degree. Specifically for any sequence $\{a(d)\geq 1\}_{d\in \mathbb{N^*}}$ one can consider $P(a(d))$, the probability that a Kostlan curve has more than $\frac{(d-1)(d-2)}{2}+1-a(d)d$ components, and there are positive constants $C$ and $D$ such that:
$$P(a(d))\leq Cd^6e^{-\frac{Dd}{a(d)}}.$$
For this random model they have also proved that the expectation of the number of components is less than $O(d)$, conjecturing that this was the right asymptotic (as confirmed by Sarnak).

Finally, for \emph{curves} there is another random model introduced by the first author in \cite{Lerario2012} that has not been discussed so far and has relation with Random Matrix Theory. Recall that once we fix three symmetric matrices $A_0, A_1, A_2$ of order $d$, we can consider the following homogeneous polynomial:
$$f(x_0, x_1, x_2)=\det(x_0A_0+x_1A_1+x_2A_2).$$
If we let the symmetric matrices to be independent random matrices in the \emph{Gaussian Orthogonal Ensemble}, $f$ becomes a random polynomial of degree $d$; its distribution is very different from the Fubini-Study ones. 
Despite this, \emph{every} real algebraic curve arises as the zero locus of one such polynomial and we can define random curves in this way; we call them \emph{random determinantal curves}. For example, using this description, maximally nested curves correspond to matrices for which there exists a positive definite linear combination $c_0A_0+c_1A_1+c_2A_2$ and since the probability of the positive definite cone becomes exponentially small in $d$, 
they are very unlikely (notice that using the above definition they also have very high energy). 

\begin{remark}
Each symmetric matrix $A$ of order $d$ also defines a quadratic form $a:\RR^d\to \RR$ by the following identity:
$$a(x)=\langle Ax, x\rangle\quad \textrm{for all $x\in \RR^d$}.$$
It is interesting that under this correspondence $A$ is in the Gaussian Orthogonal Ensemble if and only if $a$ is Kostlan distributed. Using this property the first author has computed the expectation of the sum of the Betti numbers of the intersection of one or two Kostlan distributed quadrics in $\RP^{d-1}$, showing that it behaves asymptotically as $d$ (see \cite{Lerario2012}).

Moreover there is a kind of \emph{duality} between the curve above defined and the common zero locus of $a_0, a_1, a_2$, as discussed in \cite{Complexity}. This duality allows to interchange the topology of the random intersection of three quadrics and the one of corresponding random determinantal curve, i.e. to study one using the other. 
The authors plan to discuss this ensemble and related properties in a forthcoming paper.
\end{remark}
\subsection{Betti numbers of random surfaces}Consider a random polynomial $f\in W_{3,d}$ and its zero locus $Y=Z_{S^3}(f)$ on the sphere $S^3$. Such $Y$ is a random surface and using the above results, 
combined with the computation of the expectation of its Euler characteristic, we can provide an asymptotic for \emph{each} of its Betti numbers. 
Applying Theorem 1 from \cite{Burgisser2007} with $s=1$ and $n=3$ we get:
$$\mathbb{E}\chi(Z_{\RP^3}(f))=-\frac{\sqrt{\delta}}{2}(\delta-3),$$
where $\delta$ is the parameter of the distribution, which we computed to be $\frac{d(d+2)}{3}$. Since $Y$ double covers $Z_{\RP^3}(f)$, 
its Euler characteristic is twice the one of $Z_{\RP^3}(f)$. Thus, plugging the value of the parameter in the above formula we get:
$$\mathbb{E}\chi(Y)\sim -\frac{d^3}{3\sqrt{3}}=\Theta(d^3).$$
We notice now that since $Y$ is the zero locus of a globally defined function on $S^3$, then it is orientable (a nowhere vanishing normal is given by $\textrm{grad}(f)|_{Y}$). 
We can thus apply Poincar\'e duality and get the equality $b_0(Y)=b_2(Y).$ Using the definition of the Euler characteristic we get:
$$\mathbb{E}b_0(Y)=\mathbb{E}b_2(Y)=\Theta(d^3),$$
and 
$$\mathbb{E}b_1(Y)=\mathbb{E}(b_0(Y)+b_2(Y))-\mathbb{E}\chi(Y)=\Theta(d^3).$$
We notice that in the case $f$ has \emph{even} degree then it defines a function on $\RP^3$ as well, thus $Z_{\RP^3}(f)$ is orientable and the same asymptotic holds; in the case $d$ is odd, the zero locus on the sphere is still orientable, but the one on the projective space might have nonorientable components; in particular Poincar\'e duality holds only for $\mathbb{Z}_2$ coefficients. In this case we get similar asymptotic using Betti numbers with $\mathbb{Z}_2$ coefficients.
\section{Appendix}
\subsection{Milnor's bound}
As we said in the Introduction, we only need to prove the bound for a regular $f$ (the corresponding zero locus is smooth); the set of nonregular $f$ is a proper algebraic subset of $W_{n,d}$, hence with measure (and probability) zero. In fact this proof can be modified to obtain a bound of the same order also for nonregular polynomials, using techniques from semialgebraic geometry, but this is far beyond our scope.

To start with, notice that it suffices to give a bound for the sum of the Betti numbers, since this will bound in particular $b_0(Z_{\RP^n}(f));$ 
moreover the Universal coefficients Theorem implies that the sum of the Betti numbers with \emph{integer} coefficient is bounded by the sum of the Betti numbers with $\mathbb{Z}_2$ coefficients (see \cite{Hatcher}):
$$b_0(Z_{\RP^n}(f))\leq b(Z_{\RP^n}(f), \mathbb{Z})\leq b(Z_{\RP^n}(f), \mathbb{Z}_2).$$
Now we use the so called \emph{Smith's inequality}, as presented in \cite{Wilson}. This is the statement that for a topological space $X$ with an involution $\tau:X\to X$, the inequality $b(\textrm{Fix}(\tau), \mathbb{Z}_2)\leq b(X, \mathbb{Z}_2)$ holds. In particular, since $f$ is a \emph{real} polynomial,  the complex conjugation is an involution on $Z_{\CP^n}(f)$ and Smith's inequality reads:
$$b(Z_{\RP^n}(f), \mathbb{Z}_2)\leq b(Z_{\CP^n}(f), \mathbb{Z}_2).$$
We use now the fact that $Z_{\CP^n}(f)$ is a \emph{nonsingular hypersurface} in $\CP^n$ of degree $d$. 
This property has two consequences: first, its cohomology has no torsion and $b(Z_{\CP^n}(f), \mathbb{Z}_2)=b(Z_{\CP^n}(f), \mathbb{Z})$ (Lemma 3.1 of \cite{Dimca}); 
second, such cohomology is completely determined by $d$ and $n$, and in particular:
\begin{align*}
b(Z_{\CP^n}(f), \mathbb{Z}_2) &= b(Z_{\CP^n}(f), \mathbb{Z}) \\ 
&=\frac{(d-1)^{n+1}-(-1)^{n+1}}{d}+n+(-1)^{n+1}=\Theta(d^n).
\end{align*}
The explicit formula in the above line is a consequence of the \emph{adjunction} formula, which provides the Euler characteristic of $Z_{\CP^n}(f)$ (see Exercise (3.7) in \cite{Dimca}); 
the formula for the sum of the Betti numbers immediately follows from Poincar\'e duality. Finally, since $Z_{S^n}(f)$ double covers $Z_{\RP^n}(f)$, the same bound holds for:
$$b_0(Z_{S^n}(f))\leq 2b_0(Z_{\RP^n}(f)).$$

It is interesting to notice that in the case $f$ is \emph{harmonic}, 
this bound can be obtained using Courant's Nodal Domain Theorem. 
In fact this theorem implies that the number of nodal domains of $f$ on $S^n$, i.e. $b_0(S^n\backslash Z_{S^n}(f))$, is bounded by $\Theta(d^n)$ and  the long exact sequence for the pair $(S^n, S^n\backslash Z_{S^n}(f))$ provides the bound (the reader is referred to \cite{BPR} for these effective techniques from homological algebra). 
This proof using Courant's theorem applies even when the zero set of $f$ is not regular, 
but it doesn't work for a \emph{linear combination} of harmonic functions. 

\subsection{Hyperspherical harmonics}

In this section, we review a few properties of the hyperspherical harmonics, including any relevant properties of the Gegenbauer polynomials.

First recall that the Gegenbauer polynomials $P^{(\lambda)}_{l}$ are defined using the Jacobi orthogonal polynomials $P_l^{(\alpha,\beta)}$
with $\alpha=\beta=\lambda-1/2$ \cite[Formula (4.7.1)]{Szego}. Namely,
$$P^{(\lambda)}_{l}(x) = g(\lambda,l) P_l^{(\lambda - 1/2,\lambda - 1/2)}(x),$$
where
$$ g(\lambda,l) = \frac{\Gamma(\lambda + 1/2)}{\Gamma(2 \lambda)} \frac{\Gamma(l+2 \lambda)}{\Gamma(l+\lambda+1/2)}.$$

As mentioned in Section \ref{sec:lemma2}, 
an orthonormal basis for  $H_{n,l}$ can be built up inductively.
Namely, start with an orthonormal basis for
$$\bigoplus_{m\in [0, l]}H_{n-1,m}.$$
For each such basis element, we obtain a basis element for $H_{n,l}$.

This entails that the dimension of the space $H_{n,l}$
is the same as the dimension of the space $\bigoplus_{m\in [0, l]}H_{n-1,m}$.
It is instructive to check this using the known formulas for each,
that is, to verify directly that 
$$d(n,l) = \binom{n+l}{n} - \binom{n+l-2}{n} = \sum_{m=0}^{l} \binom{n+m-1}{n-1} - \binom{n+m-3}{n-1},$$
which follows from cancellation in the sum along with the recurrence relation for binomial coefficients.

We follow \cite[Chapter 3]{Avery} in order to 
state the details resulting from the inductive construction of an orthonormal basis for $H_{n,l}$ (restricted to the sphere $S^n$) of degree $l$.
Suppose an orthonormal basis $\{Y_j(\phi)\}_{j \in I_m}$ for $H_{n-1,m}$ has already been constructed for each $m=0,1,..,l$, where $I_m$ is an index set of size 
$|I_m| = \dim H_{n-1,m} = d(n-1,m)$.
We can then obtain a spherical harmonic in $n+1$ variables (restricted to the sphere $S^n$) of degree $l$ as follows
\begin{equation}\label{spherical}
Y_{l}^m(\theta,\phi) = N_l^m (\sin\theta)^m P^{\left( \frac{n-1}{2}+m \right)}_{l-m}(\cos\theta) Y_j(\phi).
\end{equation}
In the previous formula $P^{(\lambda)}_{k}$ is the Gegenbauer polynomial, $\theta$ is the angle from the north-pole, the multi-angle $\phi\in S^{n-1}$ gives the remaining coordinates, 
and the normalization constant $N_l^m$ is given by:
\begin{equation}\label{N}
N_l^m = \sqrt{ \frac{ \Gamma( \frac{n+2m+1}{2} ) (n+2m-2)! (2l+n-1)(l-m)! }{ \sqrt{\pi} \Gamma( \frac{n+2m}{2} ) (n+2m-1) (l+m+n-2)! } } .
\end{equation}
Writing these functions for all $j \in I_m$ and for all $m = 0,1,..,d$ generates an orthonormal basis for $H_{n,l}$.

In order to produce the normalized zonal harmonic $Y_{l}(\theta)$ that was heavily utilized in Section \ref{sec:lemma1}, take $m=0$ in the formula (\ref{spherical}).
So the polynomial $Y_j(\phi) = \frac{1}{\omega_{n-1}}=\frac{\Gamma(\frac{n-1}{2})}{2\pi^{\frac{n-1}{2}}}$ is constant and: 
\begin{equation}\label{zonalGeg}
Y_{l}(\theta) = Y_{l}^0(\theta) = \frac{N_l^0}{\omega_{n-1}} P^{\left( \frac{n-1}{2} \right)}_{l}(\cos\theta),
\end{equation}
where $\omega_{n-1}$ is the $(n-1)$-dimensional volume of $S^{n-1}$ 
and
$$ N_l^0 = \sqrt{ \frac{ \Gamma( \frac{n+1}{2} ) (n-2)! (2l+n-1)l! }{ \sqrt{\pi} \Gamma( \frac{n}{2} ) (n-1) (l+n-2)! } } .$$
The Gegenbauer polynomials satisfy \cite[Formula (4.7.3)]{Szego}: 
$$P^{\left( \frac{n-1}{2}\right)}_{l}(1)=  {l + n - 2 \choose l}.$$
Together with (\ref{zonalGeg}) this implies that
\begin{equation}\label{zonalorder}
 Y_{l}(0) = \Theta(l^{(n-1)/2}).
\end{equation}

Writing this in terms of the Jacobi polynomial $P_l^{(\frac{n-2}{2},\frac{n-2}{2})}$,
\begin{equation}\label{zonal}
Y_{l}(\theta) = c(n,l) P_l^{(\frac{n-2}{2},\frac{n-2}{2})}(\cos \theta),
\end{equation}
where
\begin{equation}\label{c(n,l)}
c(n,l) = \frac{ N_l^0}{\omega_{n-1}} g \left( \frac{n-1}{2},l \right) = \Theta(l^{1/2}).
\end{equation}

Theorem 8.21.12 in \cite{Szego} states an asymptotic for the Jacobi polynomials.
In particular this applies to the polynomial $P^{(\lambda - 1/2,\lambda - 1/2)}_{l}$, 
and taking $\lambda = \frac{n-1}{2}$, we obtain the asymptotic that was used in Section \ref{sec:lemma1}.
$$P^{(\frac{n-2}{2},\frac{n-2}{2})}_{l}(\cos \theta)=\bigg(\sin \frac{\theta}{2}\cdot\cos \frac{\theta}{2}\bigg)^{\frac{2-n}{2}}\bigg\{ h(n,l)\bigg(\frac{\theta}{\sin \theta}\bigg)^{1/2}J_{\frac{2-n}{2}}(N \theta)+R_l(\theta)\bigg\},$$
where 
$$N=\frac{2l + n-1}{2},$$
\begin{equation}\label{h(n,l)} h(n,l)=\left(\frac{2l+n-1}{2} \right)^{\frac{2-n}{2}} \frac{\Gamma(l+n/2)}{l!} = \Theta(1),\end{equation}
and the error term $R_l(\theta)$ is always less than $\theta^{1/2}O(l^{-\frac{3}{2}})$ (see the remark in \cite{Szego} immediately following Theorem 8.21.12).

\subsection{Proof of the basic Claim stated in the proof of Lemma \ref{lemma2}}

Using the notation of Section \ref{sec:lemma2}, let us recall the Claim:
$$\EE \max_{S^{n-1}}  \left|F \right| \leq \sigma(m,d) \cdot d(n-1,m).$$

The proof uses a standard reproducing kernel argument.
Let $Z_m(\eta,\phi)$ be the zonal harmonic in $n$ variables with pole at the point $\phi \in S^{n-1}$.\
This is not the normalized zonal, it is the reproducing kernel as defined in \cite[p. 94]{Ax}.
Namely, $Z_m(\eta,\phi)$ acts as a reproducing kernel for the whole space of homogeneous harmonics of degree $m$.
In particular,
$$ F(\phi) = \int_{S^{n-1}} F(\eta) Z_m(\eta,\phi) d \eta .$$
The Cauchy-Schwarz inequality yields
\begin{equation}
\label{eq:CS}
  |F(\phi)|^2 \leq d(n-1,m) \int_{S^{n-1}} F(\eta)^2 d \eta ,
\end{equation}
where we have used $\int_{S^{n-1}} Z_d(\eta,\phi)^2 d \eta = d(n-1,m)$ 
\cite[p. 95, Proposition 5.27]{Ax}.

Now, taking the maximum over the sphere, then taking expectations over both sides of the inequality (\ref{eq:CS}),
we have:
\begin{align*}
\EE \max_{S^{n-1}} |F|^2 &\leq d(n-1,m) \EE \int_{S^{n-1}} F(\eta)^2 d \eta \\
&= d(n-1,m)  \sum_{j \in I_m} \EE (\hat{\xi}^j)^2 \\
&\leq d(n-1,m)^2 \sigma(m,d)^2.
\end{align*}
The inequality in the Claim now follows from the general inequality $ \EE |Y| \leq \sqrt{ \EE |Y|^2 }$ for a random variable $Y$ (it is an application of Cauchy-Schwarz again).

\end{document}